\documentclass[a4paper,11pt]{amsart}
\usepackage{graphicx}
\usepackage{mathptmx}
\usepackage{amsmath}
\usepackage{amssymb}
\usepackage{enumitem}
\usepackage{xcolor}
\usepackage{xparse}
\NewDocumentCommand{\eulerian}{omm}
 {%
  \genfrac<>{0pt}{}{#2}{#3}%
  \IfValueT{#1}{_{\!#1}}%
 }

 \newmuskip\pFqmuskip

\newcommand*\pFq[6][8]{%
  \begingroup 
  \pFqmuskip=#1mu\relax
  \mathchardef\normalcomma=\mathcode`,
  \mathcode`\,=\string"8000
  \begingroup\lccode`\~=`\,
  \lowercase{\endgroup\let~}\pFqcomma
  {}_{#2}F_{#3}{\left(\genfrac..{0pt}{}{#4}{#5}\bigg|#6\right)}%
  \endgroup
}
\newcommand{\pFqcomma}{{\normalcomma}\mskip\pFqmuskip}

\newtheorem{theorem}{Theorem}

\newtheorem{remark}[theorem]{Remark}

\begin{document}

\title[Degenerate $r$-associated Stirling numbers]{Degenerate $r$-associated Stirling numbers}

\author{Taekyun  Kim}
\address{Department of Mathematics, Kwangwoon University, Seoul 139-701, Republic of Korea}
\email{tkkim@kw.ac.kr}

\author{DAE SAN KIM}
\address{Department of Mathematics, Sogang University, Seoul 121-742, Republic of Korea}
\email{dskim@sogang.ac.kr}


\begin{abstract}
For any positive integer $r$, the $r$-associated Stirling number of the second kind $S_{2}^{(r)}(n,k)$ enumerates the number of partitions of the set $\{1,2,3,\dots,n\}$ into $k$ non-empty disjoint subsets such that each subset contains at least $r$ elements. We introduce the degenerate $r$-associated Stirling numbers of the second kind and of the first kind. They are degenerate versions of the $r$-associated Stirling numbers of the second kind and of the first kind, and reduce to the degenerate Stirling numbers of the second kind and of the first kind for $r=1$. The aim of this paper is to derive recurrence relations for both of those numbers.
\end{abstract}

\maketitle

\section{Introduction and preliminaries}
Explorations for the degenerate versions of some special numbers and polynomials have become lively interests for some mathematicians in recent years, which began from the pioneering work of Carlitz (see [1,2]). These have been done by employing various methods, such as generating functions, combinatorial methods, $p$-adic analysis, umbral calculus, operator theory, differential equations, special functions, probability theory and analytic number theory (see [5,9-13,16,17] and the references therein). \par
The Stirling number of the second kind $S_{2}(n,k)$ enumerates the number of partitions of the set $[n]=\left\{1,2,\dots,n \right\}$ into $k$ nonempty disjoint sets, while the $r$-associated Stirling number of the second kind $S_{2}^{(r)}(n,k)$ counts the number of partitions of the set $[n]$ into $k$ non-empty disjoint subsets such that each subset contains at least $r$ elements, for any positive integer $r$. The aim of this paper is to introduce the degenerate $r$-associated Stirling numbers of the second and of the first kind, and to derive their recurrence relations. They are degenerate versions of the $r$-associated Stirling numbers of the second kind and of the first kind, and reduce to the degenerate Stirling numbers of the second kind and of the first kind for $r=1$. Here we mention that the degenerate Stirling numbers of both kinds appear very frequently when one studies various degenerate versions of some special numbers and polynomials.\par
The outline of this paper is as follows. In Section 1, we recall the degenerate exponentials and the degenerate logarithms. We state the generating functions and recurrence relations for the Stirling numbers of the second kind and of the first kind.
We also remind the reader of the generating functions and recurrence relations of the degenerate Stirling numbers of the second kind and of the first kind. Finally, we recall the $r$-associated Stirling numbers of the second kind, its generating function and the recurrence relation. Section 2 is the main result of this paper. We introduce the degenerate $r$-associated Stirling numbers of the second kind and express the degenerate $r$-associated Bell polynomials in terms of the degenerate $r$-associated Stirling numbers of the second kind. In Theorem 2, we derive a recurrence relation for the degenerate Stirling numbers of the second kind. We also introduce the degenerate $r$-associated Stirling numbers of the first kind  and deduce a recurrence relation for those numbers in Theorem 4. For the rest of this section, we recall the facts that are needed throughout this paper. \par
For any nonzero $\lambda\in\mathbb{R}$, the degenerate exponentials are defined by 
\begin{equation}
e_{\lambda}^{x}(t)=(1+\lambda t)^{\frac{x}{\lambda}}=\sum_{n=0}^{\infty}\frac{(x)_{n,\lambda}}{n!}t^{n},\quad e_{\lambda}(t)=e_{\lambda}^{(1)}(t),\quad (\mathrm{see}\ [1,2]),	\label{1}
\end{equation}
where 
\begin{equation}
(x)_{0,\lambda}=1,\quad (x)_{n,\lambda}=x(x-\lambda)\cdots(x-(n-1)\lambda),\quad (n\ge 1),\quad (\mathrm{see}\ [10,11]).\label{2}
\end{equation}
Let $\log_{\lambda}t$ be the compositional inverse of $e_{\lambda}(t)$, called the degenerate logarithm, such that $\log_{\lambda}(e_{\lambda}(t))=e_{\lambda}(\log_{\lambda}(t))=t$. \par 
Then we have 
\begin{equation}
\log_{\lambda}(1+t)=\sum_{n=1}^{\infty}\frac{(1)_{n,1/\lambda}\lambda^{n-1}}{n!}t^{n},\quad (\mathrm{see}\ [9]). \label{3}	
\end{equation}
Note that $\displaystyle\lim_{\lambda\rightarrow 0}e_{\lambda}(t)=e^{t},\ \lim_{\lambda\rightarrow 0}\log_{\lambda}(1+t)=\log(1+t)$. \par 
For $k\ge 0$, the Stirling numbers of the first kind are defined by  
\begin{equation}
\frac{1}{k!}\big(\log(1+t)\big)^{k}=\sum_{n=k}^{\infty}S_{1}(n,k),\quad (k\ge 0),\quad (\mathrm{see}\ [3,13,14,17]). \label{4}	
\end{equation}
The Stirling numbers of the second kind are given by 
\begin{equation}
\frac{1}{k!}\big(e^{t}-1\big)^{k}=\sum_{n=k}^{\infty}S_{2}(n,k)\frac{t^{n}}{n!},\quad (\mathrm{see}\ [5,6,18,19,20]).\label{5}
\end{equation}
From \eqref{4} and \eqref{5}, we have 
\begin{align}
S_{1}(n+1,k)&=S_{1}(n,k-1)-nS_{1}(n,k),\label{6}\\
S_{2}(n+1,k)&=S_{2}(n,k-1)+kS_{2}(n,k),\nonumber
\end{align}
where $n,k\ge 0$ with $n\ge k$ (see [10,15,18]). \par 
Recently, the degenerate Stirling numbers of the first kind and of the second kind are respectively defined by 
\begin{equation}
\frac{1}{k!}\big(\log_{\lambda}(1+t)\big)^{k}=\sum_{n=k}^{\infty}S_{1,\lambda}(n,k)\frac{t^{n}}{n!},\label{7}
\end{equation}
and 
\begin{equation}
\frac{1}{k!}\big(e_{\lambda}(t)-1\big)^{k}=\sum_{n=k}^{\infty}S_{2,\lambda}(n,k)\frac{t^{n}}{n!},\quad (\mathrm{see}\ [9]),\label{8}	
\end{equation}
where $k$ is a non-negative integer. \par 
By \eqref{7} and \eqref{8}, we get 
\begin{align}
S_{1,\lambda}(n+1,k)&=S_{1,\lambda}(n,k-1)+(k\lambda-n)S_{1,\lambda}(n,k), \label{9}\\
S_{2,\lambda}(n+1,k)&=S_{2,\lambda}(n,k-1)+(k-n\lambda)S_{2,\lambda}(n,k), \nonumber
\end{align}
where $n,k\ge 0$ with $n\ge k$ (see [9]). \par 
Note that $\displaystyle\lim_{\lambda\rightarrow 0}S_{1,\lambda}(n,k)=S_{1}(n,k),\ \lim_{\lambda\rightarrow 0}S_{2,\lambda}(n,k)=S_{2}(n,k)\displaystyle$. \par 
The generating function of the $r$-associated Stirling numbers of the second kind is given by 
\begin{equation}
\frac{1}{k!}\bigg(e^{t}-\sum_{l=0}^{r-1}\frac{t^{l}}{l!}\bigg)^{k}=\sum_{n=rk}^{\infty}S_{2}^{(r)}(n,k)\frac{t^{n}}{n!},\quad (\mathrm{see}\ [4,6,7,8]). \label{10}
\end{equation}
Thus, by \eqref{10}, we obtain the recursion formula of $S_{2}^{(r)}(n,k),\ (n,k)\ge 0$, which is given by 
\begin{equation}
S_{2}^{(r)}(n,k)=kS_{2}^{(r)}(n-1,k)+\binom{n-1}{r-1}S_{2}^{(r)}(n-r,k-1), \label{11}	
\end{equation}
with the initial condition $S_{2}^{(r)}(n,k)=0$ if $k=0$ or $n<kr$ and $S_{2}^{(r)}(n,k)=1$ if $n=rk$. \par

\section{Degenerate $r$-associated Stirling numbers}
As degenerate versions of the $r$-associated Stirling numbers of the second kind, we consider the degenerate $r$-associated Stirling numbers of the second kind given by 
\begin{equation}
\frac{1}{k!}\bigg(e_{\lambda}(t)-\sum_{l=0}^{r-1}\frac{(1)_{l,\lambda}}{l!}t^{l}\bigg)^{k}=\sum_{n=kr}^{\infty}S_{2,\lambda}^{(r)}(n,k)\frac{t^{n}}{n!},\label{12}	
\end{equation}
where $r\in\mathbb{N}$ and $k$ is a non-negative integer. \par 
From \eqref{12}, we note that 
\begin{align}
S_{2,\lambda}^{(r)}(n,k)&=\frac{1}{k!}\sum_{\substack{l_{1}+\cdots+l_{k}=n\\ l_{i}\ge r}}\frac{n!(1)_{l_{1},\lambda} (1)_{l_{2},\lambda}\cdots (1)_{l_{k},\lambda}}{l_{1}!l_{2}!\cdots l_{k}!} \label{14} \\
&=\frac{1}{k!}\sum_{m=0}^{k}\binom{k}{m}(-1)^{m}\sum_{l_{1},l_{2},\dots,l_{m}=0}^{r-1}\frac{n!\big(\prod_{i=1}^{m}(1)_{l_{j},\lambda}\big)(k-m)_{n-l_{1}-l_{2}-\cdots-l_{m,\lambda}}}{l_{1}!l_{2}!\cdots l_{m}!(n-l_{1}-l_{2}-\cdots-l_{m})!},\nonumber
\end{align}
where $n\ge kr\ge 0,\ k,r\ge 0$. \par 
Note that $S_{2,\lambda}^{(r)}(n,k)=0$ if $n<kr$. \par 
We define the degenerate $r$-associated Bell polynomials given by 
\begin{equation}
e^{x\big(e_{\lambda}(t)-\sum_{l=0}^{r-1}\frac{(1)_{l,\lambda}}{l!}t^{l}\big)}=\sum_{n=0}^{\infty}\phi_{n,\lambda}^{(r)}(x)\frac{t^{n}}{n!}, (r\in\mathbb{N}). \label{15}	
\end{equation}
Thus, by \eqref{14} and \eqref{15}, we have 
\begin{align*}
\sum_{n=0}^{\infty}\phi_{n,\lambda}^{(r)}(x)\frac{t^{n}}{n!}&=\sum_{k=0}^{\infty}x^{k}\frac{1}{k!}\bigg(e_{\lambda}(t)-\sum_{l=0}^{r-1}\frac{(1)_{l,\lambda}}{l!}t^{l}\bigg)^{k} \nonumber \\
&=\sum_{k=0}^{\infty}x^{k}\sum_{n=kr}^{\infty}S_{2,\lambda}^{(r)}(n,k)\frac{t^{n}}{n!} \\
&=\sum_{n=0}^{\infty}\sum_{k=0}^{[\frac{n}{r}]}x^{k}S_{2,\lambda}^{(r)}(n,k)\frac{t^{n}}{n!}, 
\end{align*}
where $[x]$ denotes the greatest integer not exceeding $x$.
Therefore, we obtain the following theorem. 
\begin{theorem}
For any integers $n,r$ with $n\ge 0, r \ge 1$, we have 
\begin{displaymath}
	\phi_{n,\lambda}^{(r)}(x)=\sum_{k=0}^{[\frac{n}{r}]}x^{k}S_{2,\lambda}^{(r)}(n,k). 
\end{displaymath}	
\end{theorem}
Now, we want to find a recursion formula for the degenerate $r$-associated Stirling numbers of the second kind. \par 
Taking the derivative with respect to $t$ on both sides of \eqref{12}, we obtain
\begin{align}
\sum_{n=kr-1}^{\infty}S_{2,\lambda}^{(r)}(n+1,k)\frac{t^{n}}{n!}&=\sum_{n=kr}^{\infty}S_{2,\lambda}(n,k)\frac{t^{n-1}}{(n-1)!} \label{16} \\
&=\frac{d}{dt}\frac{1}{k!}\bigg(e_{\lambda}(t)-\sum_{l=0}^{r-1}\frac{(1)_{l,\lambda}}{l!}t^{l}\bigg)^{k}.\nonumber
\end{align}
Here we note that 
\begin{align}
	& \frac{d}{dt}\frac{1}{k!}\bigg(e_{\lambda}(t)-\sum_{l=0}^{r-1}\frac{(1)_{l,\lambda}}{l!}t^{l}\bigg)^{k}\label{17}\\
	&\quad =\frac{k}{k!}\bigg(e_{\lambda}(t)-\sum_{l=0}^{r-1}\frac{(1)_{l,\lambda}}{l!}t^{l}\bigg)^{k-1}\bigg(e_{\lambda}^{1-\lambda}(t)-\sum_{l=1}^{r-1}\frac{(1)_{l,\lambda}}{l!}lt^{l-1}\bigg) \nonumber \\
	&\quad =\frac{1}{(k-1)!}\bigg(e_{\lambda}(t)-\sum_{l=0}^{r-1}\frac{(1)_{l,\lambda}}{l!}t^{l}\bigg)^{k-1}\bigg(e_{\lambda}(t)-(1+\lambda t)\sum_{l=1}^{r-1}\frac{(1)_{l,\lambda}}{(l-1)!}t^{l-1}\bigg)\frac{1}{1+\lambda t}\nonumber.
\end{align}
Thus, by \eqref{16} and \eqref{17}, we get 
\begin{align}
&\sum_{n=kr-1}^{\infty}\bigg\{S_{2,\lambda}^{(r)}(n+1,k)+n\lambda S_{2,\lambda}^{(r)}(n,k)\bigg\}\frac{t^{n}}{n!}=(1+\lambda t)\sum_{n=kr-1}^{\infty}S_{2,\lambda}^{(r)}(n+1,k)\frac{t^{n}}{n!} \label{18} \\
&=\frac{1}{(k-1)!}\bigg(e_{\lambda}(t)-\sum_{l=0}^{r-1}\frac{(1)_{l,\lambda}}{l!}t^{l}\bigg)^{k-1}\bigg(e_{\lambda}(t)-\sum_{l=0}^{r-2}\frac{(1)_{l+1,\lambda}}{l!}t^{l}-\lambda\sum_{l=1}^{r-1}\frac{(1)_{l,\lambda}}{(l-1)!}t^{l}\bigg)\nonumber \\
&=\frac{1}{(k-1)!}\bigg(e_{\lambda}(t)-\sum_{l=0}^{r-1}\frac{(1)_{l,\lambda}}{l!}t^{l}\bigg)^{k-1}\bigg(e_{\lambda}(t)-\sum_{l=0}^{r-2}\frac{(1)_{l,\lambda}}{l!}(1-\lambda l)t^{l}-\lambda\sum_{l=1}^{r-1}\frac{(1)_{l,\lambda}}{(l-1)!}t^{l}\bigg) \nonumber\\
&=\frac{1}{(k-1)!}\bigg(e_{\lambda}(t)-\sum_{l=0}^{r-1}\frac{(1)_{l,\lambda}}{l!}t^{l}\bigg)^{k-1}\bigg(e_{\lambda}(t)-\sum_{l=0}^{r-1}\frac{(1)_{l,\lambda}}{l!}t^{l}\bigg) \nonumber \\
&+\frac{1}{(k-1)!} \bigg(e_{\lambda}(t)-\sum_{l=0}^{r-1}\frac{(1)_{l,\lambda}}{l!}t^{l}\bigg)^{k-1}\bigg(\frac{1}{(r-1)!}(1)_{r-1,\lambda}t^{r-1}\bigg) \nonumber \\
&+\frac{1}{(k-1)!}\bigg(e_{\lambda}(t)-\sum_{l=0}^{r-1}\frac{(1)_{l,\lambda}}{l!}t^{l}\bigg)^{k-1}\bigg(\lambda\sum_{l=0}^{r-2}\frac{(1)_{l,\lambda}}{l!}lt^{l}-\lambda\sum_{l=0}^{r-1}\frac{(1)_{l,\lambda}}{l!}lt^{l}\bigg) \nonumber \\
&=\frac{k}{k!}\bigg(e_{\lambda}(t)-\sum_{l=0}^{r-1}\frac{(1)_{l,\lambda}}{l!}t^{l}\bigg)^{k}+\frac{1}{(k-1)!}\bigg(e_{\lambda}(t)-\sum_{l=0}^{r-1}\frac{(1)_{l,\lambda}}{l!}t^{l}\bigg)^{k-1}\bigg(\frac{(1)_{r-1,\lambda}t^{r-1}}{(r-1)!}\bigg)\nonumber \\
&-\frac{\lambda}{(k-1)!}\bigg(e_{\lambda}(t)-\sum_{l=0}^{r-1}\frac{(1)_{l,\lambda}}{l!}t^{l}\bigg)^{k-1}\bigg(\frac{(1)_{r-1,\lambda}}{(r-1)!}(r-1)t^{r-1}\bigg) \nonumber \\
&=k\sum_{n=kr}^{\infty}S_{2,\lambda}^{(r)}(n,k)\frac{t^{n}}{n!}+\sum_{n=r(k-1)}^{\infty}S_{2,\lambda}^{(r)}(n,k-1)\frac{t^{n}}{n!} \frac{(1)_{r-1,\lambda}}{(r-1)!}t^{r-1} \nonumber \\
&-\lambda\frac{(r-1)(1)_{r-1,\lambda}}{(r-1)!}\sum_{n=r(k-1)}^{\infty}S_{2,\lambda}^{(r)}(n,k-1)\frac{t^{n}}{n!} t^{r-1}\nonumber
\end{align}
\begin{align}
&=\sum_{n=kr-1}^{\infty}\bigg\{kS_{2,\lambda}^{(r)}(n,k)+(1)_{r-1,\lambda}\binom{n}{r-1}S_{2,\lambda}^{(r)}(n-r+1,k-1)\nonumber \\
&\qquad -\lambda(r-1)(1)_{r-1,\lambda}\binom{n}{r-1}S_{2,\lambda}^{(r)}(n-r+1,k-1)\bigg\}\frac{t^{n}}{n!}.\nonumber
\end{align}
Comparing the coefficients on both sides of \eqref{18}, we have the following theorem. 
\begin{theorem}
For $n,k\ge 0$ with $n\ge kr-1$, we have 
\begin{align*}
&S_{2,\lambda}^{(r)}(n+1,k) 
=(k-n\lambda)S_{2,\lambda}^{(r)}(n,k)+(1)_{r-1,\lambda}\binom{n}{r-1}S_{2,\lambda}^{(r)}(n-r+1,k-1) \\
&\qquad\qquad\qquad-\lambda (r-1)(1)_{r-1,\lambda}\binom{n}{r-1}S_{2,\lambda}^{(r)}(n-r+1,k-1). 
\end{align*}
\end{theorem}
\begin{remark} If $r=1$ in Theorem 3, then we have 
\begin{equation}
S_{2,\lambda}^{(1)}(n+1,k)=(k-n\lambda)S_{2,\lambda}^{(1)}(n,k)+S_{2,\lambda}^{(1)}(n,k-1), \label{20}
\end{equation}
where $n,k\ge 0$ with $n\ge k-1$. 
So our result agrees with the fact in \eqref{9}, as $S_{2,\lambda}^{(1)}(n,k)=S_{2,\lambda}(n,k)$. \par
\end{remark} 
Now, we consider the degenerate $r$-associated Stirling numbers of the first kind given by  
\begin{equation}
\frac{1}{k!}\bigg(\log_{\lambda}(1+t)-\sum_{l=1}^{r-1}\frac{(1)_{r,1/\lambda}\lambda^{l-1}}{l!}t^{l}\bigg)^{k}=\sum_{n=kr}^{\infty}S_{1,\lambda}^{(r)}(n,k)\frac{t^{n}}{n!},\label{21}	
\end{equation}
where $k$ is a nonnegative integer and $r\ge 2$. \par 
Taking the derivative with respect to $t$ on both sides of \eqref{21}, we get
\begin{align}
\sum_{n=kr-1}^{\infty}S_{1,\lambda}^{(r)}(n+1,k)\frac{t^{n}}{n!} &=\sum_{n=kr}^{\infty}S_{1,\lambda}^{(r)}(n,k)\frac{t^{n}}{(n-1)!} \label{22} \\
&=\frac{d}{dt}\frac{1}{k!}\bigg(\log_{\lambda}(1+t)-\sum_{l=1}^{r-1}\frac{(1)_{l,1/\lambda}}{l!}\lambda^{l-1}t^{l}\bigg)^{k}. \nonumber	
\end{align}
Here we observe that
\begin{align}
&\frac{d}{dt}\frac{1}{k!}\bigg(\log_{\lambda}(1+t)-\sum_{l=1}^{r-1}\frac{(1)_{l,1/\lambda}}{l!}\lambda^{l-1}t^{l}\bigg)^{k} \label{23}\\
&=\frac{k}{k!}\bigg(\log_{\lambda}(1+t)-\sum_{l=1}^{r-1}\frac{(1)_{l,1/\lambda}}{l!}\lambda^{l-1}t^{l}\bigg)^{k-1}\bigg(\frac{(1+t)^{\lambda}}{1+t}-\sum_{l=1}^{r-1}\frac{(1)_{l,1/\lambda}}{(l-1)!}\lambda^{l-1}t^{l-1}\bigg).\nonumber
\end{align}
Thus, by \eqref{22} and \eqref{23}, we get 

\begin{align}
&\sum_{n=kr-1}^{\infty}S_{1,\lambda}^{(r)}(n+1,k)\frac{t^{n}}{n!}(1+t)\label{24}	\\
&=\frac{1}{(k-1)!}\bigg(\log_{\lambda}(1+t)-\sum_{l=1}^{r-1}\frac{(1)_{l,1/\lambda}}{l!}\lambda^{l-1}t^{l}\bigg)^{k-1} \bigg((1+t)^{\lambda}-\sum_{l=1}^{r-1}\frac{(1)_{l,1/\lambda}}{(l-1)!}\lambda^{l-1}t^{l-1}(1+t)\bigg) \nonumber \\
& =\frac{1}{(k-1)!}\bigg(\log_{\lambda}(1+t)-\sum_{l=1}^{r-1}\frac{(1)_{l,1/\lambda}}{l!}\lambda^{l-1}t^{l}\bigg)^{k-1}\nonumber\\
&\quad \times \bigg((1+t)^{\lambda}-\sum_{l=0}^{r-2}\frac{(1)_{l+1,1/\lambda}}{l!}\lambda^{l}t^{l}-\sum_{l=1}^{r-1}\frac{(1)_{l,1/\lambda}}{l!}l\lambda^{l-1}t^{l}\bigg) \nonumber \\
&=\frac{1}{(k-1)!}\bigg(\log_{\lambda}(1+t)-\sum_{l=1}^{r-1}\frac{(1)_{l,1/\lambda}}{l!}\lambda^{l-1}t^{l}\bigg)^{k-1}\nonumber \\
&\quad \times \bigg((1+t)^{\lambda}-\sum_{l=0}^{r-2}(1)_{l,1/\lambda}\bigg(1-\frac{l}{\lambda}\bigg)\lambda^{l}
\frac{t^{l}}{l!}-\sum_{l=1}^{r-1}\frac{(1)_{l,1/\lambda}}{l!}\lambda^{l-1}lt^{l}\bigg)\nonumber \\
&= \frac{1}{(k-1)!}\bigg(\log_{\lambda}(1+t)-\sum_{l=1}^{r-1}\frac{(1)_{l,1/\lambda}}{l!}\lambda^{l-1}t^{l}\bigg)^{k-1}\nonumber \\
&\quad \times \bigg((1+t)^{\lambda}-1-\lambda\sum_{l=1}^{r-2}\frac{\lambda^{l-1}(1)_{l,1/\lambda}}{l!}t^{l}+\sum_{l=0}^{r-2}(1)_{l,1/\lambda}\lambda^{l-1}l\frac{t^{l}}{l!}\bigg) \nonumber \\
&\quad -\frac{1}{(k-1)!}\bigg(\log_{\lambda}(1+t)-\sum_{l=1}^{r-1}\frac{(1)_{l,1/\lambda}}{l!}\lambda^{l-1}t^{l}\bigg)^{k-1}\bigg(\sum_{l=1}^{r-1}\frac{(1)_{l,1/\lambda}}{l!}\lambda^{l-1}lt^{l}\bigg) \nonumber \\
&=\frac{\lambda}{(k-1)!}\bigg(\log_{\lambda}(1+t)-\sum_{l=1}^{r-1}\frac{(1)_{l,1/\lambda}}{l!]}\lambda^{l-1}t^{l}\bigg)^{k-1}\bigg(\log_{\lambda}(1+t)-\sum_{l=1}^{r-1}\frac{(1)_{l,1/\lambda}}{l!}\lambda^{l-1}t^{l}\bigg) \nonumber \\
&\quad +\lambda^{r-1}\frac{(1)_{r-1,1/\lambda}}{(r-1)!}\frac{1}{(k-1)!}\bigg(\log_{\lambda}(1+t)-\sum_{l=1}^{r-1}\frac{(1)_{l,1/\lambda}}{l!}\lambda^{l-1}t^{l}\bigg)^{k-1} t^{r-1} \nonumber \\
&\quad -\frac{(r-1)(1)_{r-1,1/\lambda}}{(r-1)!}\lambda^{r-2}\frac{1}{(k-1)!}\bigg(\log_{\lambda}(1+t)-\sum_{l=1}^{r-1}\frac{(1)_{l,1/\lambda}}{l!}\lambda^{l-1}t^{l}\bigg)^{k-1}t^{r-1}\nonumber \\
&=k\lambda\sum_{n=kr}^{\infty}S_{1,\lambda}^{(r)}(n,k)\frac{t^{n}}{n!}+\lambda^{r-1}(1)_{r-1,1/\lambda}\sum_{n=kr-1}^{\infty}\binom{n}{r-1}S_{1,\lambda}^{(r)}(n-r+1,k-1)\frac{t^{n}}{n!} \nonumber \\
&\quad-(r-1)(1)_{r-1,1/\lambda}\lambda^{r-2}\sum_{n=kr-1}^{\infty}\binom{n}{r-1}S_{1,\lambda}^{(r)}(n-r+1,k-1)\frac{t^{n}}{n!}\nonumber \\
&=\sum_{n=kr-1}^{\infty}\bigg\{k\lambda S_{1,\lambda}^{(r)}(n,k)+\lambda^{r-1}(1)_{r-1,1/\lambda}\binom{n}{r-1}S_{1,\lambda}^{(r)}(n-r+1,k-1)\nonumber \\
&-(r-1)(1)_{r-1,1/\lambda}\lambda^{r-2}\binom{n}{r-1}S_{1,\lambda}^{(r)}(n-r+1,k-1)\bigg\}\frac{t^{n}}{n!}. \nonumber
\end{align}
On the other hand, by simple calculation, we get 
\begin{align}
&\sum_{n=kr-1}^{\infty}S_{1,\lambda}^{(r)}(n+1,k)\frac{t^{n}}{n!}(1+t)\label{25}	\\
&=\sum_{n=kr-1}^{\infty}S_{1,\lambda}^{(r)}(n+1,k)\frac{t^{n}}{n!}+\sum_{n=kr}^{\infty}S_{1,\lambda}^{(r)}(n,k)n\frac{t^{n}}{n!} \nonumber \\
&=\sum_{n=kr-1}^{\infty}\big(S_{1,\lambda}^{(r)}(n+1,k)+nS_{1,\lambda}^{(r)}(n,k)\big)\frac{t^{n}}{n!}. \nonumber
\end{align}
From \eqref{24} and \eqref{25}, we obtain the following theorem. 
\begin{theorem}
Let $r\in\mathbb{N}$ with $r\ge 2$. Then, for $n,k\ge 0$ with $n \ge kr-1$, we have 
\end{theorem}
\begin{align*}
&S_{1,\lambda}^{(r)}(n+1,k)+nS_{1,\lambda}^{(r)}(n,k) \\
&=k\lambda S_{1,\lambda}^{(r)}(n,k)+\lambda^{r-1}(1)_{r-1,1/\lambda}\binom{n}{r-1}S_{1,\lambda}^{(r)}(n-r+1,k-1)  \\
&\quad -(r-1)(1)_{r-1,1/\lambda}\lambda^{r-2}\binom{n}{r-1}S_{1,\lambda}^{(r)}(n-r+1,k-1). 
\end{align*}
\begin{remark} If $r=1$ in Theorem 4, then we have 
\begin{equation}
S_{1,\lambda}^{(1)}(n+1,k)=(k\lambda-n)S_{1,\lambda}^{(1)}(n,k)+S_{1,\lambda}^{(1)}(n,k-1), \label{20}
\end{equation}
where $n,k\ge 0$ with $n\ge k-1$. 
So our result agrees with the fact in \eqref{9}, as $S_{1,\lambda}^{(1)}(n,k)=S_{1,\lambda}(n,k)$. \par
\end{remark} 

\section{conclusion}
In recent years, studying degenerate versions of some special numbers and polynomials have drawn the attention of many mathematicians with their regained interests not only in combinatorial and arithmetical properties but also in applications to differential equations, identities of symmetry and probability theory. \par  
In this paper, we introduced the degenerate $r$-associated Stirling numbers of the second kind and of the first kind, and derived their recurrence relations. They are degenerate versions of the $r$-associated Stirling numbers of the second kind and of the first kind, and reduce to the degenerate Stirling numbers of the second kind and of the first kind for $r=1$. \par
As one of our future research projects, we would like to continue to explore degenerate versions of some special numbers and polynomials and their applications to physics, science and engineering.


\begin{thebibliography}{9}
\bibitem{1}
Carlitz, L. \emph{Degenerate Stirling, Bernoulli and Eulerian numbers.} Utilitas Math. \textbf{15} (1979), 51-88.
\bibitem{2}
Carlitz, L. \emph{A degenerate Staudt-Clausen theorem.} Arch. Math. (Basel) \textbf{7} (1956), 28-33.
\bibitem{3}
Comtet, L. \emph{Advanced combinatorics. The art of finite and infinite expansions.} Revised and enlarged edition. D. Reidel Publishing Co., Dordrecht, 1974. xi+343 pp. ISBN: 90-277-0441-4.
\bibitem{4}
Connamacher, H.; Dobrosotskaya, J. \emph{On the uniformity of the approximation for $r$-associated Stirling numbers of the second kind.} Contrib. Discrete Math. \textbf{15} (2020), no. 3, 25-42.
\bibitem{5}
Duran, U.; Acikgoz, M. \emph{On degenerate truncated special polynomials.} Mathematics, Mathematics, 2020, 8, 144.
\bibitem{6}
Fray, R. \emph{A generating function associated with the generalized Stirling numbers.} Fibonacci Quart. \textbf{5} (1967), 356-366.
\bibitem{7}
Howard, F. T. \emph{Associated Stirling numbers.} Fibonacci Quart. \textbf{18} (1980), no. 4, 303-315.
\bibitem{8}
Howard, F. T. \emph{Congruences for the Stirling numbers and associated Stirling numbers.} Acta Arith. \textbf{55} (1990), no. 1, 29-41.
\bibitem{9}
Kim, D. S.; Kim, T. \emph{A note on a new type of degenerate Bernoulli numbers.} Russ. J. Math. Phys. 27 (2020), no. 2, 227-235.
\bibitem{10}
Kim, T.; Kim, D. S. \emph{On some degenerate differential and degenerate difference operators.} Russ. J. Math. Phys. \textbf{29} (2022), no. 1, 37-46.
\bibitem{11}
Kim, T.; Kim, D. S. \emph{Some identities on truncated polynomials associated with degenerate Bell polynomials.} Russ. J. Math. Phys. \textbf{28} (2021), no. 3, 342-355.
\bibitem{12}
Kim, T.; Kim, D. S.; Jang, L.-C.; Lee, H.; Kim, H. \emph{Representations of degenerate Hermite polynomials.} Adv. in Appl. Math. 139 (2022), Paper No. 102359.
\bibitem{13}
 Kim, T.;  Kim, D. S.; Kim, H. K. \emph{Normal ordering of degenerate integral powers of number operator and its applications.} Appl. Math. Sci. Eng. \textbf{30} (2022), no. 1, 440-447.
\bibitem{14}
Kucukoglu, I. Simsek, Y. \emph{Construction and computation of unified Stirling-type numbers emerging from p-adic integrals and symmetric polynomials.} Rev. R. Acad. Cienc. Exactas Fís. Nat. Ser. A Mat. RACSAM \textbf{115} (2021), no. 4, Paper No. 167, 24 pp.
\bibitem{15}
Lehmer, D. H. \emph{Numbers associated with Stirling numbers and Number theory}  Rocky Mountain J. Math. \textbf{15} (1985), no. 2, 461-479.
\bibitem{16}
Lim, D. \emph{A note on the fully degenerate Bell polynomials of the second kind.}  Contrib. Discrete Math. \textbf{17} (2022), no. 1, 13-30.
\bibitem{17}
Park, J.-W.; Kim, B. M.; Kwon, J. \emph{Some identities of the degenerate Bernoulli polynomials of the second kind arising from $\lambda$-Sheffer sequences.} Proc. Jangjeon Math. Soc. \textbf{24} (2021), no. 3, 323-342.
\bibitem{18}
Roman, S. \emph{The umbral calculus.} Pure and Applied Mathematics, 111. Academic Press, Inc. [Harcourt Brace Jovanovich, Publishers], New York, 1984. x+193 pp.
\bibitem{19}
Simsek, B. \emph{Some identities and formulas derived from analysis of distribution functions including Bernoulli polynomials and Stirling numbers.} Filomat \textbf{34} (2020), no. 2, 521–527.
\bibitem{20}
Simsek, Y. \emph{Identities associated with generalized Stirling type numbers and Eulerian type polynomials.} Math. Comput. Appl. \textbf{18} (2013), no. 3, 251-263.
\end{thebibliography}
\end{document}